\documentclass[12pt,reqno]{amsart}
\newtheorem{lemma}{Lemma}

\newtheorem{proposition}[lemma]{Proposition}
\newtheorem{theorem}[lemma]{Theorem}
\newtheorem{remark}{Remark}[section]

\newtheorem{corollary}[lemma]{Corollary}
\newtheorem{example}{Example}[section]
\usepackage{hyperref}
\usepackage{mathtools}
\usepackage{amssymb,amsmath,amsthm}
\usepackage{epsfig,amssymb,amsmath,amsthm,graphicx,psfrag}
\usepackage[all]{xy}
\usepackage{enumerate}
\usepackage{color}
\usepackage{animate}
\usepackage{float}
\usepackage{subfig}

\usepackage{fancyhdr}
\begin{document}
	
	\title[The  group of affine transformations of homogeneous spaces with disc. isotropy]{The  group of affine transformations of homogeneous spaces with discrete isotropy }

\author{ Saldarriaga, O.$^1$, $^2$, and Fl\'orez, A.$^2$}

\subjclass[2010]{Primary: 53B05, 57S20; Secondary: 57M60, 53A15  \\ Partially Supported by CODI, Universidad de Antioquia. Project Number 2020-33713.}
\date{\today}

\begin{abstract} We present a method to compute the group of affine transformations of a 
homogeneous $G$-space under specific conditions: when the group $G$  and the homogeneous $G$-space admit linear connections so that the natural projection is affine, and with discrete  isotropy group. If $G$ admits a bi-invariant linear connection, we establish conditions under which the homogeneous space admits an invariant linear connection. As a consequence, when the isotropy group is discrete, their respective groups of affine transformations are locally isomorphic. As an application of our work, we calculate the group of the affine transformations of orientable flat affine surfaces and 3-dimensional flat affine tori.
\end{abstract}
\maketitle


\vskip5pt
\noindent
$^1$ Mathematics Department, High Point University, USA

\noindent 
$^2$ Instituto de Matem\'aticas, Universidad de Antioquia, Medell\'in-Colombia 

e-mails: osaldarr@highpoint.edu, walexander.florez@udea.edu.co 
\section{Preliminaries}

This paper deals with affine transformations of connected manifolds  endowed with a linear connection. Throughout this work manifolds are assumed to be connected, real, finite dimensional,  and without boundary (unless otherwise specified). Given a manifold $M$ endowed with a linear connection $\nabla$, we define an affine transformation of $(M,\nabla)$  as a diffeomorphism $F$ of $M$ preserving $\nabla$, that is, verifying $F_*\nabla_XY=\nabla_{F_*X}F_*Y$, for all $X,Y\in \mathfrak{X}(M),$  where $\mathfrak{X}(M)$ is the space of smooth vector fields on $M$.   The set of  diffeomorphisms preserving $\nabla$ will be  denoted by $\it{Aff}(M,\nabla)$ and it is known that under the open-compact topology and composition $\it{Aff}(M,\nabla)$ is a Lie group (see \cite{KoNo} page 229) and is called the group of affine transformations of $(M,\nabla)$.  An infinitesimal affine transformation of $(M,\nabla)$ is a smooth vector field $X$ on $M$ whose local 1-parameter groups $\phi_t$ are local affine transformations of $(M,\nabla)$.   We will denote by $\mathfrak{a}(M,\nabla)$ the real vector space of infinitesimal affine transformations of $(M,\nabla)$. 
The vector subspace $\it{aff}(M,\nabla)$ of $\mathfrak{a}(M,\nabla)$ whose elements are complete, with the usual bracket of vector fields,  is the Lie algebra of the group $\it{Aff}(M,\nabla)$  (see \cite{KoNo}).

Recall that the torsion  and curvature tensors of a connection $\nabla$ are defined by
\begin{align*} T_{\nabla}(X,Y)&=\nabla_XY-\nabla_YX-[X,Y] \\
	K_{\nabla}(X,Y)&=[\nabla_X,\nabla_Y]  - \nabla_{[X,Y]}
\end{align*}
for any $X,Y\in\mathfrak{X}(M)$. If the curvature and torsion tensors of $\nabla$ are both null, the connection is called flat affine and the pair $(M,\nabla)$ is called a flat affine manifold.  These kind of manifolds are naturally related to lagrangian foliations of symplectic manifolds (Theorem 7.8 in \cite{W}, see also \cite{F}). They are also relevant in integrable systems and mirror symmetry when their holonomy reduces to $GL_n(\mathbb{Z})$ (see \cite{KS}).

If $M=G$ is a Lie group, a linear connection on $G$ is called left invariant if every left multiplication is an affine transformation of $(G,\nabla)$, and bi-invariant if both left and right multiplications are affine transformations. Having a left invariant linear connection $\nabla^+$ on a Lie group $G$ is equivalent to have a bilinear product on $\mathfrak{g}=$Lie$(G)$ given by $X\cdot Y=(\nabla^+_{X^+}Y^+)_\epsilon$, where $X^+,Y^+$ are the left invariant vector fields on $G$ determined respectively by $X$ and $Y$ and $\epsilon$ is the group identity. When the bilinear product is given, the connection is defined by $\nabla_{X^+}Y^+=(X\cdot Y)^+$ forcing the conditions $\nabla_{fX^+}Y^+=f\nabla_{X^+}Y^+$ and $\nabla_{X^+}fY^+=X^+(f)Y^++f\nabla_{X^+}Y^+$. The connection $\nabla^+$ is torsion free if and only if  
\[ [X,Y]=X\cdot Y-Y\cdot X
\] and it is flat (i.e., with zero curvature) if and only if $[X,Y]\cdot Z=X\cdot(Y\cdot Z)-Y\cdot (X\cdot Z)$, for all $X,Y,Z\in\mathfrak{g}$. Combining these last two equations, we get that $\nabla^+$  is flat affine  if and only if the bilinear product is left symmetric, that is, 
\[ (X\cdot Y)\cdot Z-(Y\cdot X)\cdot Z=X\cdot(Y\cdot Z)-Y\cdot (X\cdot Z). \]
In this case, the pair $(G,\nabla^+)$ is called a flat affine Lie group and the algebra $(\mathfrak{g},\cdot)$  a left symmetric algebra (see \cite{Vin2}, see also \cite{JLK}). If the connection is bi-invariant and flat affine, the algebra $(\mathfrak{g},\cdot)$ is associative (see \cite{Med}).

If $p:\widehat{M}\rightarrow M$ is the universal covering map of a real $n$-dimensional flat affine manifold $(M,\nabla)$. The pullback $\widehat{\nabla}$  of $\nabla$ by $p$ is a flat affine structure on $\widehat{M}$ and $p$ is an affine map. Moreover, the group  of deck transformations, which in this case is isomorphic to $\pi_1(M)$, acts on $\widehat{M}$ by affine transformations.  There exists an affine immersion $D:(\widehat{M},\widehat{\nabla})\longrightarrow (\mathbb{R}^n,\nabla^0)$, with $\nabla^0$ the usual linear connection on $\mathbb{R}^n$, and a  group homomorphism  $A:\it{Aff}(\widehat{M},\widehat{\nabla}) \longrightarrow \it{Aff}(\mathbb{R}^n,\nabla^0)$ so that the following diagram commutes 
 	\[ \xymatrix{ \widehat{M} \ar[d]_{F} \ar[r]^{D} &\mathbb{R}^n\ar[d]^{A(F)}\\
 		\widehat{M} \ar[r]^{D} &\mathbb{R}^n.} \]
 The map $D$ is called  a developing map and it was introduced by Ehresmann (see \cite{E}). In particular for every $\gamma\in\pi_1(M)$, we have $D\circ \gamma= \mathsf{h}(\gamma)\circ D$ with $\mathsf{h}(\gamma):=A(\gamma)$. The map $\mathsf{ h}$ is also a group homomorphism called the holonomy representation of $(M,\nabla)$.   
 
 It is known that $(G,\nabla^+)$ is a flat affine Lie group if and only if  there exists an affine  \'etale representation $\rho:\widehat{G}\rightarrow\it{Aff}(\mathbb{R}^n)$ (see \cite{JLK} and \cite{Med}). This means that the respective action of $\widehat{G}$ on $\mathbb{R}^n$ leaves an open orbit $\mathcal{O}$ with discrete isotropy. The open orbit turns out to be the image of a developing map $D:\widehat{G}\rightarrow \mathbb{R}^n$ and $D$ is a covering map of $\mathcal{O}$.  

Finally, given a Lie group $G$, a homogeneous $ G $-space $ M $ is a manifold admitting a transitive action $\tau:G\times M\rightarrow M$. The isotropy group $ H $ at any point $ p\in M $ is a closed subgroup of $G$, hence it is a Lie subgroup and the set of left cosets $G/H$ admits a unique structure of manifold so that the projection $\pi:G\rightarrow G/H$ is a smooth map. Moreover, the manifolds $G/H$ and $M$ are diffeomorphic under the identification $gH\mapsto \tau(g,p)$. A linear connection $\nabla$ on $M$ is called invariant if $\tau_g$ is affine relative to $\nabla$, for every $g\in G$, where $\tau_g:M\rightarrow M$ is the map defined by $\tau_g(m)=\tau(g,m)$. Invariant linear connections on a homogeneous $G$-space $M$ are in one-to-one correspondence with linear maps $ \mathcal{L}:\mathfrak{g}\rightarrow End(\mathfrak{g/h}) $ satisfying
 \begin{align} \label{Eq:invariantconnectionsonhomogeneousspaces1}
		& \mathcal{L}_X(Y+\mathfrak{h})=[X,Y]+\mathfrak{h}\qquad\qquad\qquad\  \text{ for all } X\in\mathfrak{h}, Y\in\mathfrak{g} \\  \label{Eq:invariantconnectionsonhomogeneousspaces2}
	& \mathcal{L}_{Ad_h X}(Ad_hY+\mathfrak{h})=\overline{Ad}_h(\mathcal{L}_X(Y+\mathfrak{h})), \text{ for all } h\in H \text{ and }X,Y\in\mathfrak{g}.
\end{align}
where $\mathfrak{g}=$Lie$(G)$, $\mathfrak{h}=$Lie$(H)$, $Ad$ denotes the adjoint map and $\overline{Ad}_h:\mathfrak{g/h}\mapsto\mathfrak{g/h}$ the map defined by $\overline{Ad}_h(X+\mathfrak{h})={Ad}_h(X)+\mathfrak{h}$ (see \cite{Vin}, see also \cite{El}).

 A  homogeneous space $G/H$ is called reductive if there is a linear subspace $\mathfrak{m}$ of  $\mathfrak{g}$ so that $\mathfrak{g}=\mathfrak{h}\oplus \mathfrak{m}$ and $Ad_h(\mathfrak{m})\subseteq\mathfrak{m}$, for every $h\in H$. This implies that $ad_h(\mathfrak{m})\subseteq \mathfrak{m}$, for all $h\in H$. When $H$ is connected, the converse is also true. In this special case, invariant linear connections on $G/H$ are in one-to-one correspondence with bilinear maps $\theta:\mathfrak{m}\times\mathfrak{m}\rightarrow\mathfrak{m}$ so that   $Ad_h|_{\mathfrak{m}}\subseteq Aut(\mathfrak{m},\theta)$, for any $h\in H$ (see \cite{Nom}, see also \cite{El}).
 
This work is organized as follows. In Section \ref{S:affineconnectionsongspacesfromconnectionsong}, we study homogeneous spaces $G/H$ where $G$ is endowed with a left invariant linear connection $\nabla^+$ and give conditions on $\mathfrak{h}=$Lie$(H)$ so that $G/H$ admits an invariant connection determined by $\nabla^+$. We also give conditions on the connection so that the natural projection is an affine map. The reductive case is treated separately in Section \ref{S:connectionsonreductivespaces}. We use Section \ref{S:transformationsonhomogeneousspacesfromtransformationsong} to study the group of affine transformations of a homogeneous space $G/H$ endowed with an invariant connection whose natural projection map is an affine map.  We prove that if $H$ is discrete, the group of affine transformations of $G/H$ is locally isomorphic to the group of affine transformations of $G$ commuting with the action of $H$ on $G$. We devote Section \ref{S:affinetransformationsofflataffinesurfaces} to calculate the group of affine transformations of the flat affine two dimensional oriented surfaces and of the 3-dimensional flat affine tori. Finally, we finish the work with some further results.  

\section{Invariant connections on homogeneous $G$-spaces from left invariant connections on $G$} \label{S:affineconnectionsongspacesfromconnectionsong}

In this section we give necessary conditions for a left invariant connection on a Lie group, to determine an invariant connection on a homogeneous $G$-manifold, see Theorem \ref{T:conexioninvariantesobreespacioshomogeneos}. We also give a version of Proposition 3.3 in \cite{Pos} (see also Lemma 1.8.24 in \cite{Wol}) for homogeneous spaces (see Theorem \ref{T:piisaffineifnablaisbiinvariant}). 

Let $ G$  be a Lie group endowed with a left invariant connection $\nabla^+$, consider a closed Lie subgroup $H$ of $G$ and let $G/H $ be the respective homogeneous manifold. Set  $ \mathfrak{g}=Lie(G)=T_\epsilon G $, $ \mathfrak{h}=Lie(H)=T_\epsilon H $ and denote by $ X^+$ (respectively $ X^- $) the left (respectively right) invariant vector field on $G$ determined by $ X $. Now consider the bilinear product on $\mathfrak{g}$ given by $X*Y=(\nabla^+_{X^-}Y^-)_\epsilon$ and the  linear map  $ \mathcal{L}:\mathfrak{g}\rightarrow End(\mathfrak{g/h}) $ assigning to every $Y\in\mathfrak{g}$ the map $\mathcal{L}_Y$ defined by 
\begin{equation} \label{Eq:mapdefiningtheprojectedconnection} \mathcal{L}_Y(X+\mathfrak{h})= X* Y+\mathfrak{h}=(\nabla^+_{X^-}Y^-)_\epsilon+\mathfrak{h}. \end{equation}
Notice that for any $Y\in\mathfrak{g}$, the map $\mathcal{L}_Y$ is well defined if and only if $X* Y\in\mathfrak{h}$, for all $X\in\mathfrak{h}$ if and only if $\mathfrak{h}$ is a right ideal of the algebra $(\mathfrak{g},*)$. In this case, according to Vinberg (see \cite{Vin}), the map $\mathcal{L}$ determines an invariant connection on $ G/H $ if and only if it verifies Equations \eqref{Eq:invariantconnectionsonhomogeneousspaces1} and \eqref{Eq:invariantconnectionsonhomogeneousspaces2}. The first equation is  equivalent to have that $Y* X-[X,Y]\in \mathfrak{h}$ for all $X\in\mathfrak{h}$.  If $ R_h\in \it{Aff}(G,\nabla^+) $ for all $ h\in H $, where $R_h$ denotes right multiplication by $h$ on $G$, that is, $R_h(g)=gh$. One can easily verify that $Ad_h(Y*X)=Ad_hY*Ad_hX$, for any $X,Y\in\mathfrak{g}$ and $h\in H$. Hence, Equation \eqref{Eq:invariantconnectionsonhomogeneousspaces2} holds. The following remark is very relevant in the sequel (see \cite{El} page 215)

\begin{remark} \label{R:flowofprojection} A vector field $ X\in\mathfrak{X}(G) $ whose flow $ \phi_t^X $ commutes with right multiplication $R_h$, for all $ h\in H $, determines a vector field $\overline{X} $ on $ G/H $ defined by
	\begin{equation} \label{Eq:flujoenelcociente} \overline{X}_{\sigma H}=\left.\dfrac{d}{dt} \right|_{t=0} \phi_t^X(\sigma)H\qquad\text{for all }\sigma\in G.  \end{equation}
	As the vector field is invariant under $R_h$ if and only if its flow commutes with $R_h$, it follows that the vector field $\overline{X}$ is well defined on $G/H$. It is also easy to see that $\overline{X}=\pi_*(X) $, where $\pi:G\rightarrow G/H$ is the natural projection. In particular, for $X\in\mathfrak{g}$, the vector field $\pi_*(X^-)$ is well defined in $G/H$ and will be denoted by  $X^*$.
\end{remark}


Recall that a connection $\overline{\nabla}$ on the homogeneous space $G/H$ is called invariant if the left action, denoted by $\tau$, of $G$ on $G/H$ preserves the connection. That is, $(\tau_g)_*\overline{\nabla}_XY=\overline{\nabla}_{(\tau_g)_*X}(\tau_g)_*Y$ for  all $X,Y\in \mathfrak{X}(G/H)$ and all $g\in G$, where $\tau_g$ is  given by $\tau_g(g'H)=(gg')H$. In fact, $\tau_g$ is the map so that the following diagram commutes 
\[ \xymatrix{ G \ar[d]_{\pi} \ar[r]^{L_g} &G\ar[d]^{\pi}\\
	G/H \ar[r]^{\tau_g} &G/H,} \]
where $L_g$ denotes left multiplication by $g$ on $G$. 
 
If $R_h\in\it{Aff}(G,\nabla^+)$ for any $h\in H$,  the vector field $ \nabla_{X^-}^+Y^- $ is $ R_h $-invariant. Hence the vector field $ \pi_*(\nabla_{X^-}^+Y^-) $ is well defined on $G/H$. Now consider the operator $\overline{\nabla}$ given by 
\begin{equation}\label{Eq:connectioninducedfromG}
	\overline{\nabla}_{X^*}Y^*:=\pi_*(\nabla_{X^-}^+Y^-)
\end{equation}
It is easy to verify that $(L_g)_*(X^-)=(Ad_gX)^-$ and that $(Ad_gX)^*=(\tau_g)_*(X^*)$, hence we have  
\begin{equation} \label{Eq:invarianceofoverlinenabla} \begin{array}{ll}
		(\tau_g)_*\overline{\nabla}_{X^*}Y^*  & =(\tau_g)_*\circ\pi_*(\nabla_{X^-}^+Y^-)=\pi_*\circ(L_g)_*(\nabla_{X^-}^+Y^-)\\
		& =\pi_*(\nabla_{(L_g)_*X^-}^+(L_g)_*Y^-)=\overline{\nabla}_{\pi_*(Ad_gX)^-}\pi_*(Ad_gY)^-\\
		& =\overline{\nabla}_{(Ad_gX)^*}(Ad_gY)^* =\overline{\nabla}_{(\tau_g)_*X^*}(\tau_g)_*Y^*.
\end{array}\end{equation}
That is, $\overline{\nabla}$ is $\tau_g$-invariant. Moreover, as the map $\phi:\mathfrak{g}/\mathfrak{h}\rightarrow T_{\epsilon H}(G/H)$ is defined by $\phi(X+\mathfrak{h})=\pi_{*,\epsilon}(X)$, and by noticing that $ \pi_{*,\epsilon}(\nabla_{X^-}^+Y^-)_{\epsilon }=\phi(\mathcal{L}_Y(X+\mathfrak{h})) $, we get that $\overline{\nabla}$ is the connection on $ G/H $ determined by the map defined in \eqref{Eq:mapdefiningtheprojectedconnection}.  From all of the previous we have the following.

\begin{theorem} \label{T:conexioninvariantesobreespacioshomogeneos} Given a  Lie group $G$ endowed with a  left invariant linear connection $\nabla^+$ and $H$  a Lie subgroup  of $G$ verifying 
	\begin{enumerate}[a.] \item $\mathfrak{h}=$Lie$(H)$ is a right ideal of the algebra $(\mathfrak{g},*)$ where $\mathfrak{g}=\text{Lie}(G)$ and $*$ as defined in \eqref{Eq:mapdefiningtheprojectedconnection}.
		\item $Y* X-[X,Y]\in \mathfrak{h}$ for all $X\in\mathfrak{h}$ and 
		\item $R_h\in\it{Aff}(G,\nabla^+)$, for all $h\in H$, \end{enumerate}
	then the connection $\overline{\nabla}$  defined by Equation \eqref{Eq:connectioninducedfromG} is an invariant connection on $G/H$. 
\end{theorem}

\begin{example} Consider the group of affine motions of the line $G=\it{Aff}(\mathbb{R})\cong\mathbb{R}^{*}\times \mathbb{R}$ with the product $(x,y)(x',y')=(xx',xy'+y)$,  where $\mathbb{R}^*$ is the multiplicative group of nonzero real numbers. Its Lie algebra is $\mathfrak{g}=\it{aff}(\mathbb{R})=\mathbb{R}e_1\oplus\mathbb{R}e_2$ with bracket $[e_1,e_2]=e_2$. The 1-dimensional linear subspaces of $\mathfrak{g}$ are given by $\mathfrak{h}=\mathbb{R}f$ with $f=e_1+\beta e_2$, $\beta\in\mathbb{R}$ or $f=e_2$. Next we exhibit all left invariant linear connections $\nabla^+$ on $G$ that determine invariant connections $\overline{\nabla}$ on $G/H$, with $H$ the Lie subgroup of $G$ of Lie algebra $\mathfrak{h}$.
	
	If $\mathfrak{h}=\mathbb{R}(e_1+\beta e_2)$ with $\beta$ fixed, then $H$  is given by $H=\{(a,\beta(a-1))\mid a\ne0\}$. A calculation shows that all left invariant connections on $G$ so that conditions a., b., and c. are verified, are  given by the bilinear product on $\it{aff}(\mathbb{R})$ determined by the table on the left 
	\begin{equation} \label{Eq:primerejemploenaffR} \begin{array}{c|c|c}  & f &e_2\\\hline  f& \alpha f& e_2\\\hline e_2&0&0	\end{array}\qquad\qquad\qquad \begin{array}{c|c|c} * & f &e_2\\\hline f& \alpha f& 0\\\hline e_2&e_2&0	\end{array} \end{equation} 
	the table on the right determines the product $*$, where $f=e_1+\beta e_2$. Notice that all these connections are flat affine. One can verify that $f^*=\pi_*(f^-)=u\frac{\partial}{\partial u}$ and $e_2^*=\frac{\partial}{\partial u}$, where $(u)$ is a system of local  coordinates on $G/H$. Hence from the product $*$ above, the connection $\overline{\nabla}$ on $G/H$ is determined by $\overline{\nabla}_{\frac{\partial}{\partial u}}\frac{\partial}{\partial u}=0$, that is, $\overline{\nabla}$ is the usual connection on $G/H$.

	Notice also that, as $\overline{\nabla}$ is determined by the map $\mathcal{L}$, it must also satisfy that $\overline{\nabla}_{e_2^*}f^*=e_2^*$, that is, $\overline{\nabla}_{\frac{\partial}{\partial u}}u\frac{\partial}{\partial u}=\frac{\partial}{\partial u}$ and this is immediately  verified using the properties of linear connections.
	
	Now, if $\mathfrak{h}=\mathbb{R}e_2$, the group   $H$ is given by $H=\{(1,a)\mid a\in\mathbb{R}\}$ and all left invariant connections on $G$ verifying the conditions of Theorem \ref{T:conexioninvariantesobreespacioshomogeneos} are determined by the bilinear product obtained from the following table displaying the values on the linear basis $(e_1,e_2)$ of $\mathfrak{g}$
	\begin{equation} \label{Eq:segundoejemploenaffR} \begin{array}{c|c|c}  & e_1 &e_2\\\hline  e_1& \alpha e_1+\lambda e_2& \gamma e_2\\\hline e_2&(\alpha-\gamma) e_2&0	\end{array}\qquad\qquad\qquad \begin{array}{c|c|c} * & e_1 &e_2\\\hline e_1& \alpha e_1+\lambda e_2& (\gamma-1)e_2\\\hline e_2&(\alpha-\gamma+1)e_2&0	\end{array}   \end{equation}
	The table on the right hand side gives the product $*$. A simple calculation shows that $e_1^*=u\frac{\partial}{\partial u}$ and $e_2^*=0$, hence the connection $\overline{\nabla}^{\alpha}$ is determined by $\overline{\nabla}^\alpha_{u\frac{\partial}{\partial u}}u\frac{\partial}{\partial u}=\alpha u\frac{\partial}{\partial u}$. The other condition given by the map $\mathcal{L}$ is $\overline{\nabla}_{e_1^*}e_2^*=(\gamma-1)e_2^*$, i.e.,  $\overline{\nabla}_{u\frac{\partial}{\partial u}}0=0$ which is obviously satisfied. 
\end{example}

\begin{remark} Under the conditions of the previous theorem, if $\nabla^+$ is torsion free, Conditions a. and b. of the theorem agree. Hence the connection $\overline{\nabla}$ exists if conditions a. and c. are satisfied.  
\end{remark}

Moreover, if $\nabla^+$ is bi-invariant, one can show that the natural projection $\pi:G\rightarrow G/H$ is an affine map, more precisely we get the following. 

\begin{theorem} \label{T:piisaffineifnablaisbiinvariant} If $G$ is a Lie group endowed with a bi-invariant linear connection $\nabla^+$ and $H$ is a Lie subgroup so that  conditions a. and b. of the previous theorem hold, then  the connection $\overline{\nabla}$  defined by Equation \eqref{Eq:connectioninducedfromG} is an invariant connection on $G/H$ so that the natural projection $\pi:\left(G^{op},\nabla^+\right)\rightarrow \left(G/H,\overline{\nabla}\right)$ is an affine map, where $G^{op}$ denotes the opposite Lie group of $G$.
\end{theorem}
\begin{proof} As $\nabla^+$ is bi-invariant, right multiplications $R_g$ belong to $\it{Aff}(G,\nabla^+)$, for every $g\in G$, hence Condition c. of Theorem \ref{T:conexioninvariantesobreespacioshomogeneos} holds. Then from the previous theorem it follows that there is an invariant connection $\overline{\nabla}$ on $G/H$ determined by Equation \eqref{Eq:connectioninducedfromG}. Moreover, since right invariant vector fields on $G$ are left invariant vector fields on $G^{op}$, the formula $\nabla^+_{X^-}Y^-$ determines a left invariant connection on $G^{op}$. Hence Equation \eqref{Eq:connectioninducedfromG} means  that $\pi:(G^{op},\nabla^+)\rightarrow \left(G/H,\overline{\nabla}\right)$ is an affine map.
\end{proof}
A similar result was proved in \cite{AbHa} (Proposition 5.7) in the context of reductive homogeneous spaces.
\begin{example} Consider the group $G=\it{Aff}(\mathbb{R})$ endowed with the connection $\nabla^{+,\alpha}$ determined by the bilinear product on the table in the left hand side of \ref{Eq:primerejemploenaffR}, for $\alpha$ fixed. One can verify that the connection $\nabla^{+,1}$ is bi-invariant, hence the projection is an affine map. In particular, notice that the connection $\overline{\nabla}$ verifies all relations of the table on the right hand side  in \ref{Eq:primerejemploenaffR}.
	
	On the other hand, a left invariant connection $\nabla^{+,\alpha,\gamma,\lambda}$ determined by a  product on the left hand side table in \eqref{Eq:segundoejemploenaffR} is bi-invariant if and only if $\lambda=0$. Hence the map $\pi:\left(G^{op},\nabla^{+,\alpha,\gamma,0} \right)\rightarrow \left(G/H,\overline{\nabla}^{\alpha}\right)$  is an affine map. However, it can be observed that $\pi$ is affine for any $\alpha$, $\gamma$, and $\lambda$. 
\end{example}

\begin{corollary} Under the conditions of Theorem \ref{T:conexioninvariantesobreespacioshomogeneos}, if $H$ is a normal subgroup, $\overline{\nabla}$ is left invariant.
\end{corollary}
\begin{proof} An easy verification shows that the left multiplication map $L_{gH}$ by $gH$ on $G/H$ coincides with the map $\tau_g$, for any $g\in G$. 
\end{proof}

\begin{example} The subgroup $H=\{(1,a)\mid a\in\mathbb{R}\}$ is a normal subgroup of $G=\it{Aff}(\mathbb{R})$, hence the connections $\overline{\nabla}$ on $G/H$ determined by the table on the left hand side in \ref{Eq:segundoejemploenaffR} are left invariant.
\end{example}

\section{Reductive case} \label{S:connectionsonreductivespaces} In this section we exhibit a different proof for Proposition 5.7 in \cite{AbHa}. Recall that a homogeneous space $M\cong G/H$ is called reductive if there is a decomposition of $\mathfrak{g}$ as a direct sum of vector spaces $\mathfrak{g}=\mathfrak{h}\oplus\mathfrak{m}$ so that $Ad_h(\mathfrak{m})\subseteq\mathfrak{m}$ for all $h\in H$. As $\pi_{*,\epsilon}X=0$ whenever $X\in\mathfrak{h}$, and $\mathfrak{m}$ can be identified with $\mathfrak{g/h}$, we have an isomorphism  $\phi:\mathfrak{m}\rightarrow T_{\epsilon H}(G/H)$ defined by $\phi(X)=\pi_{*,\epsilon}X$. Moreover, for $X=X_\mathfrak{m}+X_\mathfrak{h}\in\mathfrak{g}$, it holds that $\phi^{-1}\circ\pi_{*,\epsilon}X=X_\mathfrak{m}$. Consider the product on $\mathfrak{m}$ defined by 
\begin{equation} \label{Eq:productoreductivo} X\cdot Y=(X*Y)_{\mathfrak{m}}, \end{equation}
where $X*Y=(\nabla^+_{X^-}Y^-)_\epsilon$. Since $X\cdot Y=\phi^{-1}\circ\pi_{*,\epsilon}(X*Y)$, Equation \eqref{Eq:productoreductivo} defines a bilinear product on $\mathfrak{m}$.  Under these terms we have.

\begin{proposition} \label{T:conexioninvariantesobreespacioshomogeneosreductivos} Let $G$ be a Lie group endowed with a left invariant connection $\nabla^+$, let $H$ be a closed Lie subgroup of $G$ and  $G/H$ the respective homogeneous manifold. If $G/H$ is  reductive,  $Ad_h(\mathfrak{h})\subseteq\mathfrak{h}$,  and $R_h\in \it{Aff}(G,\nabla^+)$, for every $h\in H$,   then the product defined  in Equation \eqref{Eq:productoreductivo} determines an invariant connection on $G/H$. The connection is given by the formula
	\[ \overline{\nabla}_{X^*}Y^*=\pi_*(\nabla^+_{X^-}Y^-), \quad\text{for all}\quad X,Y\in\mathfrak{m}.  \]
Moreover, if $\nabla^+$ is bi-invariant, the natural projection $\pi:\left(G^{op},\nabla^+\right)\rightarrow \left(G/H,\overline{\nabla}\right)$ is an affine map.
\end{proposition}
\begin{proof} First notice that since $G/H$ is reductive we get that $Ad_h((X*Y)_\mathfrak{m})\in\mathfrak{m}$, and since $H$ is a Lie subgroup, we have that $Ad_h((X*Y)_\mathfrak{h})\in\mathfrak{h}$, thus	
\[ Ad_h(X*Y)=Ad_h((X*Y)_\mathfrak{m}+(X*Y)_\mathfrak{h})=\underbrace{Ad_h((X*Y)_\mathfrak{m})}_{\in\mathfrak{m}}+\underbrace{Ad_h((X*Y)_\mathfrak{h})}_{\in\mathfrak{h}}.\]
It follows that $Ad_h(X*Y)_\mathfrak{m}=Ad_h((X*Y)_\mathfrak{m})$. Now, as $L_h,R_h\in\it{Aff}(G,\nabla^+)$, for any $h\in H$, we also get   that $Ad_h(X*Y)=Ad_h(X)*Ad_h(Y)$. Hence 
	\begin{align*}Ad_h(X\cdot Y)&=Ad_h((X*Y)_\mathfrak{m})=(Ad_h(X*Y))_\mathfrak{m}=(Ad_hX*Ad_hY)_\mathfrak{m}=
	Ad_hX\cdot Ad_hY \end{align*}
%
%
%
%
%
%
%
%
	it follows that $Ad_H|_\mathfrak{m}\subseteq \it{Aut}(\mathfrak{m},\cdot)$. And, from Theorem 8.1 in \cite{Nom}, this implies that the product in Equation \eqref{Eq:productoreductivo} determines an invariant connection on $G/H$. From Equation \eqref{Eq:invarianceofoverlinenabla} and by noticing that 
	\[ \pi_{*,\epsilon}(\nabla^+_{X^-}Y^-)_\epsilon=\pi_{*,\epsilon}(X*Y)=\pi_{*,\epsilon}((X*Y)_\mathfrak{m}+(X*Y)_\mathfrak{h})
	=\pi_{*,\epsilon}((X*Y)_\mathfrak{m})=\pi_{*,\epsilon}(X\cdot Y), \]
	we conclude that the connection determined by  Equation \eqref{Eq:productoreductivo} is $\overline{\nabla}$. Finally, by mimicking the proof of Theorem \ref{T:piisaffineifnablaisbiinvariant}, we get that if $\nabla^+$ is bi-invariant and	$\pi$ is  an affine map.
\end{proof}

\section{Group of affine transformations in homogeneous spaces} \label{S:transformationsonhomogeneousspacesfromtransformationsong}

In this section we study the relationship between the group of affine transformations of a Lie group and a  homogeneous $G$-space  respectively endowed with linear connections $\nabla$ and $\overline{\nabla}$ so that the isotropy is discrete and $\pi$ is an affine map. This relationship between the groups of affine transformations, also applies to Lie groups and homogeneous spaces with discrete isotropy verifying the conditions of Section \ref{S:affineconnectionsongspacesfromconnectionsong}. We start with the following remark which will be useful in the section. 
 


\begin{remark}\label{R:diffeoprojection} Given a smooth map $\phi:G\rightarrow G$ commuting with $R_h$, for all $h\in H$, one can verify that the map $\overline{\phi}:G/H\rightarrow G/H$ given by $\overline{\phi}(gH)=\phi(g)H$ is a well defined smooth map. It is also easy to check  that $\overline{\phi}$ is one-to-one (respectively onto) whenever ${\phi}$ is one-to-one (respectively onto). 
	Notice that $\overline{\phi}$ is the unique map so that the following diagram commutes
	\begin{equation} \label{Eq:affinetransformationbelow} \xymatrix{ G \ar[d]_{\pi} \ar[r]^{\phi} &G\ar[d]^{\pi}\\
			G/H \ar[r]^{\overline{\phi}} &G/H.} \end{equation}
It is easy to check that a diffeomorphism $\phi$ defines a unique diffeomorphism $\overline{\phi}$ modulo $\it{Aut}(\pi)$. More specifically, $\overline{\phi}=\overline{\psi}$ if and only if $\phi=\psi\circ f$ for some $f\in\it{Aut}(\pi)$
\end{remark}

\begin{proposition} \label{P:Generalcase} If $G$ is a Lie group, $G/H$ a homogeneous space both  endowed with linear  connections  $\nabla$ and  $\overline{\nabla}$ and    $\pi:(G,\nabla)\rightarrow (G/H,\overline{\nabla})$ is an affine map, then we have 
	\begin{equation}\label{Eq:weakerversionoftheorem7} \left\{\overline{\phi}\mid \text{there exists }\phi\in \it{Aff}(G,\nabla) \text{ so that } \pi\circ\overline{\phi}=\phi\circ \pi \right\}\subseteq \it{Aff}(G/H,\overline{\nabla}) \end{equation}
\end{proposition}
\begin{proof} First notice that for $X\in\mathfrak{g}$ we have that
	\begin{align} \label{Eq:overlinephiisaffine} \notag \overline{\phi}_*\overline{\nabla}_{X^*}Y^*&=\overline{\phi}_*\pi_*\nabla_{X^-}Y^-=\pi_*\phi_*\nabla_{X^-}Y^-
		=\pi_*\nabla_{\phi_*X^-}\phi_*Y^-\\&=\overline{\nabla}_{\pi_*\phi_*X^-}\pi_*\phi_*Y^-
		=\overline{\nabla}_{\overline{\phi}_*\pi_*X^-}\overline{\phi}_*\pi_*Y^-=\overline{\nabla}_{\overline{\phi}_*X^*}\overline{\phi}_*Y^*. 
	\end{align}
Now, since the space $\{X^*\mid X\in\mathfrak{g}\}$ generates $\mathfrak{X}(G/H)$, it follows that $\overline\phi$ is an affine diffeomorphism.
\end{proof}

From now on, the set $\{ \phi \in\it{Aff}(G,\nabla)\mid \phi \text{ commutes with }R_h, \text{ for all }h\in H \}$ will be denoted by $\it{Aff}(G,\nabla)^\pi$. Similarly, $\it{aff}(G,\nabla)^\pi=\{X\in\it{aff}(G,\nabla)\mid (R_h)_*(X)=X\}.$ 

\begin{remark} It is easy to verify that the map $\phi\rightarrow\overline{\phi}$ defines a group homomorphism between $\it{Aff}(G,\nabla)^\pi$ and $\it{Aff}(G/H,\overline{\nabla})$ with kernel given by $\it{Aut}(\pi)$.
\end{remark}

The other inclusion does not hold in general. However, if $H$ is discrete, we can prove the equality in the infinitesimal case. The following lemma is key on its proof.  

\begin{lemma} \label{L:sameprojectionsamefield} Let $X\in\mathfrak{X}(G)$ be $R_h-$invariant, for every $h\in H$, if $\pi_*(X)\equiv0$ and $H$ is discrete, then $X\equiv0$
\end{lemma} 
\begin{proof} If $\phi_t$ is the flow of $X$, by Remark \ref{R:flowofprojection}, there is a well defined vector $\pi_*(X)$ in $\mathfrak{X}(G/H)$ with flow given by $\overline{\phi_t}(gH)=\phi_t(g)H$, for all $g\in G$, and since $\pi_*(X)\equiv0$, we have that $\overline{\phi_t}\equiv Id_{G/H}$. It follows that
	\[ \phi_t(g)H=gH \]
Hence, for every $g\in G$, $g^{-1}\phi_t(g)=h_g(t)\in H$. So, we have that $\dfrac{d}{dt}\bigg|_{t=0}g^{-1}\phi_t(g)\in\mathfrak{h}=\{0\}$. This means that $h_g(t)=h_g$ does not depend on $t$. Thus, for every $g\in G$, we have
\[ X_g=\frac{d}{dt}\bigg|_{t=0}\phi_t(g)=\frac{d}{dt}\bigg|_{t=0}gh_g=0 \]
\end{proof}

If the homogeneous space $G/H$ is reductive with $\mathfrak{g}=\mathfrak{h}\oplus\mathfrak{m}$, given a vector field $\overline{X}\in\mathfrak{X}(G/H)$, there always exists a natural lift $L\left(\overline{X}\right)\in\mathfrak{X}(G)$ so that $\pi_*(L\left(\overline{X}\right))=\overline{X}$ (see \cite{AbHa} pg 247). For this we use the horizontal distribution $\mathcal{D}$ given by $\mathcal{D}_g=(L_g)_{*,\epsilon}\mathfrak{m}$. Since $\mathfrak{m}$ is $Ad_h$ invariant, it follows that this distribution is $R_h$ invariant for every $h\in H$. Also notice that if $H$ is a discrete subgroup of $G$, the homogeneous space $G/H$ is clearly reductive, so we have the following.

\begin{proposition}\label{afinenafin} Under the same conditions of Proposition \ref{P:Generalcase}, if $H$ is a discrete subgroup of $G$, then we have
	\[ \mathfrak{a}\left(G/H,\overline{\nabla}\right) = \left\{\overline{X}\in\mathfrak{X}(G/H)\mid {L\left(\overline{X}\right)}\in\mathfrak{a}(G,\nabla) \right\}, \]
	where $L\left(\overline{X}\right)$ is the lift of $\overline{X}$.
\end{proposition}
\begin{proof} For the first inclusion, let $\overline{X}\in \mathfrak{a}(G/H,\overline{\nabla})$ and $X=L\left(\overline{X}\right)$ its lift with flow $\phi_t$. Hence the flow $\overline{\phi}_t$ of $\overline{X}$ verifies that $\overline{\phi}_t(gH)=\phi_t(g)H$, for all $g\in G$. Moreover, for any pair of right invariant vector fields $Y^-$ and $Z^-$, we have that
\begin{align*} \pi_*(\phi_t)_*\nabla_{Y^-}Z^-&=\overline{\phi}_t\pi_*\nabla_{Y^-}Z^- =\overline{\phi}_t\overline{\nabla}_{\pi_*Y^-}\pi_*Z^-
=\overline{\nabla}_{\overline{\phi}_t\pi_*Y^-}\overline{\phi}_t\pi_*Z^-
\\ &=\overline{\nabla}_{\pi_*(\phi_t)_*Y^-}\pi_*(\phi_t)_*Z^-=\pi_* \nabla_{(\phi_t)_*Y^-}(\phi_t)_*Z^-
\end{align*}
Hence, from Lemma \ref{L:sameprojectionsamefield} we get that $(\phi_t)_*\nabla_{Y^-}Z^-=\nabla_{(\phi_t)_*Y^-}(\phi_t)_*Z^-$, that is $X\in \mathfrak{a}(G,\nabla)$.

For the other inclusion, take $\overline{X}\in \mathfrak{X}(G/H)$ so that $X=L\left(\overline{X}\right)\in \mathfrak{a}(G,\nabla)$ and let $\overline{\phi}_t$ and $\phi_t$ be their respective flows. Since $\pi\circ\phi=\overline{\phi}\circ\pi$, we have that
\begin{align*} (\overline{\phi}_t)_*\overline{\nabla}_{Y^*}Z^*&=(\overline{\phi}_t)_*\overline{\nabla}_{\pi_*Y^-}\pi_*Z^-=
(\overline{\phi}_t)_*\pi_*{\nabla}_{Y^-}Z^-=\pi_*({\phi}_t)_*{\nabla}_{Y^-}Z^-
\\&=\pi_*{\nabla}_{({\phi}_t)_*Y^-}({\phi}_t)_*Z^-={\overline{\nabla}}_{\pi_*({\phi}_t)_*Y^-}(\pi_*{\phi}_t)_*Z^- 
\\&={\overline{\nabla}}_{(\overline{\phi}_t)_*\pi_*Y^-}(\overline{\phi}_t)_*\pi_*Z^-
={\overline{\nabla}}_{(\overline{\phi}_t)_*Y^*}(\overline{\phi}_t)_*Z^*
\end{align*}
Therefore $\overline{X}\in\mathfrak{a}(G/H,\overline{\nabla})$
\end{proof}

\begin{remark} Notice that when $H$ is a discrete subgroup of $G$,  if a vector field $X\in\mathfrak{X}(G)$ is projectable to a vector field $\overline{X}=\pi_*(X)$,  from Lemma \ref{L:sameprojectionsamefield}, we have that $X=L(\overline{X})$. Hence $X$ is invariant under $R_h$ for every $h\in H$. Therefore, from Remark \ref{R:flowofprojection} we get that a vector field is projectable if and only if it is $R_h-$invariant for every $h\in H$. 
\end{remark}

Recall that a covering $\left(\overline{M},\pi,M\right)$ is regular if automorphisms of the covering act transitively on its fibers. It is known that when $H$ is discrete, as the action of $H$ on $G$ is smooth, the map $\pi:G\rightarrow G/H$ determines a regular covering (see Example 7.1 in \cite{Lim}, see also \cite{Pos} page 34). Hence we have the following.

\begin{theorem} \label{C:affinetransformationgrouponhomogeneousspaceswithdiscretequotient} If $G$ is a Lie group, $H$ a discrete subgroup of $G$, $G$ and $G/H$ are respectively endowed with linear connections $\nabla$ and $\overline{\nabla}$, and $\pi:G\rightarrow G/H$ is affine, then we have
	\[ \it{aff}\left(G/H,\overline{\nabla}\right) = \left\{\overline{X}\in\mathfrak{X}(G/H)\mid {L\left(\overline{X}\right)}\in\it{aff}(G,\nabla) \right\}, \]
	where $L\left(\overline{X}\right)$ is the lift of $\overline{X}$.
\end{theorem}
\begin{proof} If $\overline{X}\in\mathfrak{X}(G/H)$ is so that $L\left(\overline{X}\right)\in\it{aff}(G,\nabla)$ with corresponding flows $\overline{\phi}_t$ and $\phi_t$, we get that $\phi_t\in\it{Aff}(G,\nabla)$ for every $t$, and by Remark \ref{R:diffeoprojection} and Proposition \ref{P:Generalcase}, $\overline{\phi}_t\in\it{Aff}(G/H,\overline{\nabla})$  for every $t$, thus  $\overline{X}\in\it{aff}(G/H,\overline{\nabla})$.

For the other inclusion, let $\overline{X}\in\it{aff}\left(G/H, \overline{\nabla}\right)$ with flow $\overline{\phi}_t$ and let $\phi_t$ be the flow of $L\left(\overline{X}\right)$. We claim that $\phi_t$ is a diffeomorphism for every $t$. Let $D$ be a fundamental domain for the right action of $H$ over $G$ and let $D'\subset D$ be a subset containing exactly one representative for every class $gH$, with $g\in G$. Notice that $\phi_t(D')$ also contains exactly one representative for every class because if $\phi_t(g_1)H=\phi_t(g_2)H$, with $g_1,g_2\in D'$, it follows that $\overline{\phi}_t(g_1H)=\overline{\phi}_t(g_2H)$. But, since $\overline{\phi}_t$ is injective, we have that $g_1H=g_2H$, therefore $g_1=g_2$. Thus  $\pi_{/D'}$ and $\pi_{/\phi_t(D')}$ are both bijective, so $\phi_t{_{/D'}}$ is also a bijection. Hence, using also the fact that the covering is regular, we get that $\phi_t$ is bijective and the claim follows. Finally, from Proposition \ref{afinenafin} we conclude that $L\left(\overline{X}\right)\in\it{aff}(G,\nabla)$ 
\end{proof}

\begin{remark} Under the conditions of the previous theorem, it follows from  Lemma \ref{L:sameprojectionsamefield} that the map $X\rightarrow\pi_*(X)$ gives a natural isomorphism between  $\it{aff}(G,\nabla)^\pi$ and  $\it{aff}(G/H,\overline{\nabla})$ with inverse $\overline{X}\mapsto L\left(\overline{X}\right)$.  Hence we have the following
\end{remark}
\begin{corollary} \label{C:affinetransformationsofG/HareaffinetransformationsofGcommutingwithRh} Under the same assumptions of the previous theorem, we get that the group $\it{Aff}(G,\nabla)^\pi/\it{Aut}(\pi)$ is locally isomorphic to  $\it{Aff}(G/H,\overline{\nabla})$.  
\end{corollary}
\begin{proof} Since the map $\phi\rightarrow\overline{\phi}$ is a homomorphism from $\it{Aff}(G,\nabla)^\pi$ to $\it{Aff}(G/H,\overline{\nabla})$ with kernel $\it{Aut}(\pi)$, we get that $\it{Aff}(G,\nabla^+)^\pi/\it{Aut}(\pi)$ is isomorphic to a subgroup of $\it{Aff}(G/H,\overline{\nabla})$. The result follows by noticing that $\it{Aff}(G,\nabla)^\pi/Aut(\pi)$ have the same dimension as $\it{Aff}^\pi(G/H,\overline{\nabla})$ as shown in the previous remark. 
\end{proof}

Using Theorem \ref{T:piisaffineifnablaisbiinvariant}, If $\nabla$ is bi-invariant, then some of the conditions in the preceding theorem are redundant. In particular, the existence of a linear connection on $G/H$ and the fact that the projection $\pi$ is affine, are direct consequences of the bi-invariance of $\nabla$. Therefore we get

\begin{theorem} \label{T:affinetransformationgrouponhomogeneousspaceswithdiscretequotient} If $G$ is a Lie group endowed with a bi-invariant connection $\nabla^+$ and $H$ is a discrete subgroup of $G$, then 
\[ \it{aff}\left(G/H,\overline{\nabla}\right) = \left\{X\in\mathfrak{X}(G/H)\mid L(X)\in\it{aff}(G,\nabla^+) \right\}, \] 
where $\overline{\nabla}$ is the invariant connection on $G/H$ guaranteed by Theorem \ref{T:piisaffineifnablaisbiinvariant}. Moreover, the group $\it{Aff}(G,\nabla)^\pi/\it{Aut}(\pi)$  is locally isomorphic to  $\it{Aff}\left(G/H,\overline{\nabla}\right) $.
\end{theorem}
\begin{proof} Just notice that since $\mathfrak{h}=$Lie$(H)$ is trivial, conditions a. and b. of Theorem \ref{T:conexioninvariantesobreespacioshomogeneos} hold, hence by Theorem \ref{T:piisaffineifnablaisbiinvariant}, there exists a connection $\overline{\nabla}$ on $G/H$ so that $\pi$ is an affine map. Hence the first assertion follows from Theorem \ref{C:affinetransformationgrouponhomogeneousspaceswithdiscretequotient} and the previous corollary. The previous corollary also implies the second statement.
\end{proof}

\section{Group of affine transformations of  flat affine surfaces} \label{S:affinetransformationsofflataffinesurfaces}

As an application of Theorem \ref{C:affinetransformationgrouponhomogeneousspaceswithdiscretequotient}, in the following examples we calculate  the  groups of affine transformations of all flat affine orientable surfaces. We need the following remark

\begin{remark} \label{R:groupsofaffinetransformationsoforbits} From Lemmas 4 and 5 in \cite{SF}, we get  that the groups of classical affine transformations preserving the upper half plane $\mathcal{H}_2$, the quadrant $\mathcal{C}_2=\{(x,y)\mid x,y>0\}$ and the punctured plane are respectively given by 
	$ G_1=\left\{ T:\mathbb{R}^2\rightarrow\mathbb{R}^2\mid T(x,y)=(ax+by+c,dy),\ d>0 \right\}$, $ G_2=\left\{ T:\mathbb{R}^2\rightarrow\mathbb{R}^2\mid T(x,y)=(ax,by)\text{ with }a,b>0 \right\}$ and $G_3=GL_2(\mathbb{R}).$
\end{remark}

\begin{example} \label{E:affinetwotorii} (Group of affine transformations of the flat affine tori) It is well known that, up to isomorphism, there are six flat affine homogeneous structures on the 2-torus $\mathbb{T}^2$ (see \cite{Nag}, \cite{AF} and \cite{Ben}). These are determined respectively by  \'etale affine representations of $\mathbb{R}^2 $, $ \rho_i:\mathbb{R}^2\longrightarrow \it{Aff}(\mathbb{R}^2) $, $i=1,\dots,6$ defined by
	\[ \begin{array}{lll}
		\rho_1(a,b)=\begin{bmatrix}
			1 & 0 & a \\
			0 & 1 & b \\
		\end{bmatrix},\quad & \rho_2(a,b)=\begin{bmatrix}
			1 & b & a+\frac{1}{2}b^2 \\
			0 & 1 & b \\
		\end{bmatrix},\quad &
		\rho_3(a,b)=\begin{bmatrix}
			1 & 0 & a \\
			0 & e^b & 0 \\
		\end{bmatrix},\\\\ \rho_4(a,b)=\begin{bmatrix}
			e^a & be^a & 0 \\
			0 & e^a & 0 \\
		\end{bmatrix},\quad &
		\rho_5(a,b)=\begin{bmatrix}
			e^a & 0 & 0 \\
			0 & e^b & 0 \\
		\end{bmatrix},\quad & \rho_6(a,b)=e^a\begin{bmatrix}
			\ \ \cos b & \sin b & 0 \\
			-\sin b & \cos b & 0 \\
		\end{bmatrix}
	\end{array} \]
where each of the above matrices represent the affine transformation of $\mathbb{R}^2$ whose linear part is given by the first two columns  and the last column corresponds to the translation part.  Each of these representations determines an action of $\mathbb{Z}^2$ on the open orbit  whose quotient is a  torus.  The actions are as follows 
	
	\noindent  Case 1.  $\mathbb{Z}^2$ acting on the plane by $\xi_1((m,n),(x, y))=(x+m,y+n)$.
	
	\noindent	Case 2.  The quotient obtained by the action   $\xi_2((m, n),(x, y))=(x+ny+ m+n^2/2,y+n)$ of $\mathbb{Z}^2$ on the plane. This torus is isomorphic to the Kuiper torus.

	\noindent Case 3. The action of $\mathbb{Z}^2$ on the upper half plane $\mathcal{H}_2$ given by $\xi_3((m, n),(x, y))=(x+m,e^ny).$

	\noindent Case 4. $\mathbb{Z}^2$ acting  on  $\mathcal{H}_2$ by  $\xi_4((m, n),(x, y))=e^m(x+ny,y)$.

	\noindent Case 5. The action of $\mathbb{Z}^2$ on  the quadrant $\mathcal{C}_2$  given  by $\xi_5((m, n),(x, y))=(e^mx,e^ny).$

	\noindent Case 6. The quotient of the punctured plane $M=\mathbb{R}^2\setminus\{(0,0)\}$ by the action of $\mathbb{Z}$ given by $\xi_6(m,(x, y))=e^m(x,y)$. This is known as the  Hopf torus. 
	
	The following picture shows each of these tori
	
	\begin{figure}[H]
		\centering
		\subfloat[\tiny{Usual torus (Case 1)}]{
			\label{f:xxx}
			\includegraphics[width=0.27\textwidth]{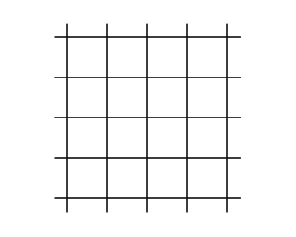}}
		\subfloat[\tiny{Kuiper torus (case 2)}]{
			\label{f:xxx}
			\includegraphics[width=0.27\textwidth]{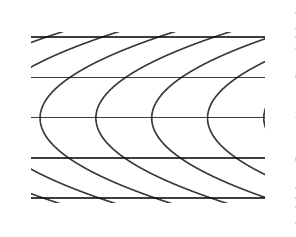}}
		\subfloat[\tiny{Case 3}]{
			\label{f:xxx}
			\includegraphics[width=0.27\textwidth]{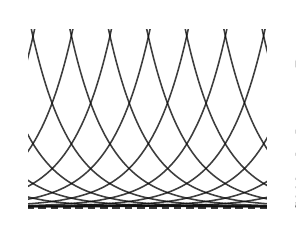}}\\ \end{figure} \begin{figure}[H] 
		\subfloat[\tiny{Case 4}]{
			\label{f:xxx}
			\includegraphics[width=0.26\textwidth]{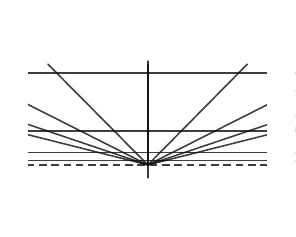}}
		\subfloat[\tiny{Case 5}]{
			\label{f:xxx}
			\includegraphics[width=0.26\textwidth]{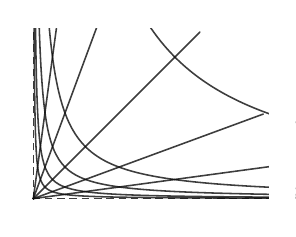}}
		\subfloat[\tiny{Hopf torus (Case 6)}]{
			\label{f:xxx}
			\includegraphics[width=0.26\textwidth]{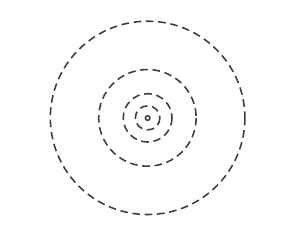}}
		\caption{Regular covering of the 2-torus}
	\end{figure}	

\vskip-10pt
For $i=1,\dots,6$, a developing map $D_i$ induces a group structure on $D_i(\mathbb{R}^2)$. In particular, for $i=6$, taking $D_6(x,y)=e^x(cos(2\pi y),\sin(2\pi y))$, one obtains the product of nonzero complex numbers. For $i\leq5$ that, $D_i(\mathbb{R}^2)/D_i(\mathbb{Z}^2)$ is a homogeneous space isomorphic to the torus $\mathbb{T}_i$. Notice that $D_6(\mathbb{Z}^2)=\{ (e^m,0)\mid m\in\mathbb{Z} \}$, hence the torus $\mathbb{T}_6$ is isomorphic to $\mathbb{C}^*/\mathbb{Z}$ where $\mathbb{C}^*$ is the group of nonzero complex numbers and $\mathbb{Z}$ is identified with the subgroup $\{(e^m,0)\mid m\in\mathbb{Z}\}$ of $\mathbb{C}^*$. By letting $\nabla^i$ be the left invariant flat affine connection on $\mathbb{R}^2$ determined by $\rho_i$.  By Theorem \ref{C:affinetransformationgrouponhomogeneousspaceswithdiscretequotient}, there exists an invariant connection $\overline{\nabla}^i$ on the torus $\mathbb{T}_i$, $i=1,\dots, 6$ and, by Corollary \ref{C:affinetransformationsofG/HareaffinetransformationsofGcommutingwithRh}. the group of affine transformations is locally isomorphic to the group of classical affine transformations preserving the orbit and commuting with the corresponding action. Using  the groups of transformations preserving the upper half plane, the quadrant and the punctured plane given in the previous remark and also noticing that the $\it{Aut}(\pi)\cong\mathbb{Z}^2$ where $\pi:\mathbb{R}^2\rightarrow\mathbb{T}_i$ is the natural projection, for $i=1,\dots,5$ and $\it{Aut}(\pi)\cong\mathbb{Z}$ when $i=6$, we get \\
	
	$\it{Aff}(\mathbb{T}_1,\overline{\nabla}^1)$  is isomorphic to the usual torus, the group of translations of $\mathbb{R}^2/\mathbb{Z}^2$.
	
	$ \it{Aff}(\mathbb{T}_2,\overline{\nabla}^2)\cong \{F:\mathbb{R}^2\rightarrow\mathbb{R}^2 \mid F(x,y)=(x+ay+b,y+a)\}/\mathbb{Z}^2.  $
	
	$\it{Aff}(\mathbb{T}_3,\overline{\nabla}^3)\cong \{F:\mathbb{R}^2\rightarrow\mathbb{R}^2 \mid F(x,y)=(x+a,by),\  b>0\}/\mathbb{Z}^2\cong (\mathbb{R}^{>0},\cdot)/\mathbb{Z}\times (\mathbb{R},+)/\mathbb{Z} $   where $\mathbb{R}^{>0}$ 
	
	is the set of positive real numbers.
	
	$\it{Aff}(\mathbb{T}_4,\overline{\nabla}^4)\cong\{F:\mathbb{R}^2\rightarrow\mathbb{R}^2 \mid F(x,y)=(ax+by,ay),\  a>0\}/\mathbb{Z}^2\cong(\mathbb{R}^{>0},\cdot)/\mathbb{Z}\times (\mathbb{R},+)/\mathbb{Z}. $  
	
	$\it{Aff}(\mathbb{T}_5,\overline{\nabla}^5)\cong\{F:\mathbb{R}^2\rightarrow\mathbb{R}^2 \mid F(x,y)=(ax,by),\  a,b>0\}/\mathbb{Z}^2\cong (\mathbb{R}^{>0},\cdot)/\mathbb{Z}\times (\mathbb{R}^{>0},\cdot)/\mathbb{Z}.$
	
	$ \it{Aff}(\mathbb{T}_6,\overline{\nabla}^6)\cong  GL_2(\mathbb{R})/\mathbb{Z}.$\\

The first five cases are abelian and $2-$dimensional and the last one is non-abelian and 4-dimensional, as it is mentioned without a proof in \cite{Nag}. 

\noindent	
 Also notice that each group $\it{Aff}(\mathbb{T}^2,\overline{\nabla}^i)$, $i=1,\dots,5$ is compact isomorphic to the torus itself. 
\end{example}

\begin{example} (Groups of affine transformations of flat affine cylinders) In all the following cases, the quotient of the given space by the action of $\mathbb{Z}$,  determines an affine cylinder
	\begin{enumerate}[(i)]
		\item The plane with the action defined by $\psi_1(m,x,y)=(x,y+m)$.
		\item The action on the plane defined by $\psi_2(m,x,y)=(x+my+m^2/2,y+m)$.
		\item The upper half plane and the action given by $\psi_3(m,x,y)=(x+m,y)$.
		\item The action on the upper half plane defined as $\psi_4(m,x,y)=(x,e^my)$.
		\item The upper half plane with the action $\psi_5(m,x,y)=(x+my,y)$.
		\item The action on the upper half plane given by $\psi_6(m,x,y)=(e^mx,e^my)$.
		\item The cuadrant with the action $\psi_7(m,x,y)=(x,e^my)$.
		\item The action $\psi_8(m,x,y)=(x,y+m)$ on the plane with the connection $\nabla^6$. This  cylinder is affinely isomorphic to the punctured plane.
	\end{enumerate}
Cylinders corresponding to the cases (i) through (vii) are diffeomorphic to the homogeneous space $D(\mathbb{R}^2)/D(\mathbb{Z})$ where $D$ is a developing map and $\mathbb{Z}$ is identified either with the discrete subgroup $\{(m,0)\mid m\in\mathbb{Z}\}$ or $\{(0,m)\mid m\in\mathbb{Z}\}$. The last one is diffeomorphic to $\mathbb{R}^2/\mathbb{Z}$. Hence using Remark \ref{R:groupsofaffinetransformationsoforbits}, Theorem \ref{C:affinetransformationgrouponhomogeneousspaceswithdiscretequotient}, and Corollary \ref{C:affinetransformationsofG/HareaffinetransformationsofGcommutingwithRh}  we get that the respective groups of affine transformations  are given by \\
	
	$\it{Aff}(\mathcal{C}_1,\overline{\nabla}'^1)\cong\{ F:\mathbb{R}^2\rightarrow\mathbb{R}^2\mid F(x,y)=(ax+b,cx+y+d),\ a\ne0 \}/\mathbb{Z}.$ 
	
	$\it{Aff}(\mathcal{C}_2,\overline{\nabla}'^2)\cong\{ F:\mathbb{R}^2\rightarrow\mathbb{R}^2\mid F(x,y)=(x+ay+b,y+a) \}/\mathbb{Z}.$
	
	$\it{Aff}(\mathcal{C}_3,\overline{\nabla}'^3)\cong\{ F:\mathbb{R}^2\rightarrow\mathbb{R}^2\mid F(x,y)=(x+ay+b,cy),\ c>0 \}/\mathbb{Z}.$
	
	$\it{Aff}(\mathcal{C}_4,\overline{\nabla}'^4)\cong\{ F:\mathbb{R}^2\rightarrow\mathbb{R}^2\mid F(x,y)=(ax+b,cy),\ a\ne0,\ c>0 \}/\mathbb{Z}.$
	
	$\it{Aff}(\mathcal{C}_5,\overline{\nabla}'^5)\cong\{ F:\mathbb{R}^2\rightarrow\mathbb{R}^2\mid F(x,y)=(ax+by+c,ay),\ a>0 \}/\mathbb{Z}.$
	
	$\it{Aff}(\mathcal{C}_6,\overline{\nabla}'^6)\cong\{ F:\mathbb{R}^2\rightarrow\mathbb{R}^2\mid F(x,y)=(ax+by,cy),\ a\ne0,\ c>0 \}/\mathbb{Z}.$
	
	$\it{Aff}(\mathcal{C}_7,\overline{\nabla}'^7)\cong\{ F:\mathbb{R}^2\rightarrow\mathbb{R}^2\mid F(x,y)=(ax,by),\ a,b>0 \}/\mathbb{Z}.$
	
	$\it{Aff}(\mathcal{C}_8,\overline{\nabla}'^8)\cong GL_2(\mathbb{R})$
	
	\noindent
	where $\overline{\nabla}'^i$ is the connection on the cylinder determined from the corresponding connection on $\mathbb{R}^2$.
\end{example}

\subsection{Three dimensional Tori} We use the  classification of the three dimensional flat affine structures on $\mathbb{R}^3$, due Remm-Goze (see \cite{ReGo}), to find 15 non-isomorphic flat affine structures on the three dimensional torus. Then we apply Theorem \ref{C:affinetransformationgrouponhomogeneousspaceswithdiscretequotient} to find the algebra of affine transformations and Corollary \ref{C:affinetransformationsofG/HareaffinetransformationsofGcommutingwithRh}  to find a group  locally isomorphic to the group of affine transformations. Remm-Goze's classification yields the following connections on $\mathbb{R}^3$ determined by the given developing maps, which we slightly modified to get appropriate orbits of $\mathbb{R}^3$ to use Lemma 4. in \cite{SF}. We present them in the following table where we exhibit their corresponding developing maps, a basis for $\it{aff}(\mathbb{R},\nabla_i)$, the action determining the torus $\mathbb{T}_i$, a basis for $\it{aff}^\pi(\mathbb{R},\nabla_i)$, i.e, vector fields in $\it{aff}(\mathbb{R},\nabla_i)$ commuting with the action, and a matrix group locally isomorphic to $\it{Aff}(\mathbb{T}_i,\overline{\nabla}_i)$, for $i=1,\dots,15$ (the last column of the $3\times4$ matrices corresponds to the translation part)

\noindent 
\begin{tabular}{|c|c|}\hline 
$D_1(a,b,c)=(e^a,e^{a+b},e^{a+c}),\ \mathcal{O}_1=\left(\mathbb{R}^{>0}\right)^3$ & $D_2(a,b,c)=(e^a\cos b,e^{a}\sin b,e^{a}(e^c-\cos b))$
\\
$\left\{x\frac{\partial }{\partial x},y\frac{\partial }{\partial y},z\frac{\partial }{\partial z}\right\}$ & $\left\{x\frac{\partial }{\partial x},y\frac{\partial }{\partial x},x\frac{\partial }{\partial y},y\frac{\partial }{\partial y},x\frac{\partial }{\partial z},y\frac{\partial }{\partial z},z\frac{\partial }{\partial z},\frac{\partial }{\partial z}\right\}$  \\ 
$(m,n,p)\cdot(x,y,z)=(e^mx,e^{m+n}y,e^{m+p}z)$ & $(m,n,p)\cdot(x,y,z)=(e^mx,e^my,e^m(e^p-1)x+e^{m+p}z)$ \\

$\left\{x\frac{\partial }{\partial x},y\frac{\partial }{\partial y},z\frac{\partial }{\partial z}\right\}$ & $\left\{ x\frac{\partial }{\partial x}+z\frac{\partial }{\partial z},x\frac{\partial }{\partial y},y\frac{\partial }{\partial y} \right\}$ \\ 
$\it{Aff}(\mathbb{T}_1,\overline{\nabla}_1)_0\cong\left\{\left[\begin{smallmatrix}
	a&0&0\\0&b&0\\0&0&c
\end{smallmatrix}\right]\bigg|a,b,c\ne0\right\}/\mathbb{Z}^3$ &\  
$\it{Aff}(\mathbb{T}_2,\overline{\nabla}_2)_0=\left\{\left[\begin{smallmatrix}
	a&0&0\\b&c&0\\0&0&a
\end{smallmatrix}\right]\ \bigg|\ a,c\ne0 \right\}/\mathbb{Z}^2$\\ \hline 
$D_3(a,b,c)=(ae^c,e^{b+c},e^c)$ & $D_4(a,b,c)=\left(\left(a+\frac{b^2}{2}\right)e^c,be^{c},e^c\right)$  \\
 $\left\{x\frac{\partial }{\partial x},y\frac{\partial }{\partial x},z\frac{\partial }{\partial x},\frac{\partial }{\partial x},y\frac{\partial }{\partial y},z\frac{\partial }{\partial z}\right\}$ & $\left\{x\frac{\partial }{\partial x},y\frac{\partial }{\partial x},z\frac{\partial }{\partial x},\frac{\partial }{\partial x},x\frac{\partial }{\partial y},y\frac{\partial }{\partial y},z\frac{\partial }{\partial y},\frac{\partial }{\partial y},z\frac{\partial }{\partial z}\right\}$  \\
$(m,n,p)\cdot(x,y,z)=(e^px+me^pz,e^{n+p}y,e^pz)$ & $(m,n,p)\cdot(x,y,z)=$\\ 
\ &$(e^px+ne^py+(m+\frac{n^2}{2})e^{p}z,e^py+ne^pz,e^pz)$ \\
 $\left\{ x\frac{\partial }{\partial x}+z\frac{\partial }{\partial z},z\frac{\partial }{\partial x},y\frac{\partial }{\partial y} \right\}$ & $\left\{ x\frac{\partial }{\partial x}+y\frac{\partial }{\partial y}+z\frac{\partial }{\partial z},y\frac{\partial }{\partial x}+z\frac{\partial }{\partial y},z\frac{\partial }{\partial x}, \right\}$ \\ 
$\it{Aff}(\mathbb{T}_3,\overline{\nabla}_3)_0\cong \left\{\left[\begin{smallmatrix}
		a&0&b\\0&c&0\\0&0&a
	\end{smallmatrix}\right]\ \big|\ a,c\ne0  \right\}/\mathbb{Z}^3$ & $\it{Aff}(\mathbb{T}_4,\overline{\nabla}_4)_0\cong \left\{\left[\begin{smallmatrix}
		a&b&c\\0&a&b\\0&0&a
	\end{smallmatrix}\right]\ \bigg|\ a\ne0  \right\}/\mathbb{Z}^3$ \\ \hline 
 $D_5(a,b,c)=(ae^c,be^{c},e^c)$ & $D_6(a,b,c)=(e^a,e^{b},c)$ \\
 $\left\{x\frac{\partial }{\partial x},y\frac{\partial }{\partial x},z\frac{\partial }{\partial x},\frac{\partial }{\partial x},x\frac{\partial }{\partial y},y\frac{\partial }{\partial y},z\frac{\partial }{\partial y},\frac{\partial }{\partial y},z\frac{\partial }{\partial z}\right\}$ & $\left\{x\frac{\partial }{\partial x},y\frac{\partial }{\partial y},x\frac{\partial }{\partial z},y\frac{\partial }{\partial z},z\frac{\partial }{\partial z},\frac{\partial }{\partial z}\right\}$ \\
 $(m,n,p)\cdot(x,y,z)=(e^px+me^pz,e^py+ne^pz,e^pz)$ & $(m,n,p)\cdot(x,y,z)=(e^mx,e^ny,z+p)$  \\ 
 $\left\{ x\frac{\partial }{\partial x}+y\frac{\partial }{\partial y}+z\frac{\partial }{\partial z},z\frac{\partial }{\partial x},z\frac{\partial }{\partial y}, \right\}$ & $\left\{ x\frac{\partial }{\partial x},y\frac{\partial }{\partial y},\frac{\partial }{\partial z} \right\}$ \\
 $\it{Aff}(\mathbb{T}_5,\overline{\nabla}_5)_0\cong \left\{\left[\begin{smallmatrix}
 		a&0&c\\0&a&b\\0&0&a
 	\end{smallmatrix}\right]\ \bigg|\ a\ne0  \right\}/\mathbb{Z}^3.$ & $\it{Aff}(\mathbb{T}_6,\overline{\nabla}_6)_0\cong \left\{\left[\begin{smallmatrix}
 		a&0&0&0\\0&b&0&0\\0&0&1&c
 	\end{smallmatrix}\right]\ \bigg|\ a,b\ne0  \right\}/\mathbb{Z}^3.$ \\ \hline 
\end{tabular}
\newpage\noindent   \begin{tabular}{|c|c|}\hline
 	$D_7(a,b,c)=(e^a\cos b,e^a\sin b,c)$ & $D_8(a,b,c)=(e^a, b,c)$ \\
 	$\left\{x\frac{\partial }{\partial x},y\frac{\partial }{\partial x},x\frac{\partial }{\partial y},y\frac{\partial }{\partial y},x\frac{\partial }{\partial z},y\frac{\partial }{\partial z},z\frac{\partial }{\partial z},\frac{\partial }{\partial z}\right\}$ & $\left\{x\frac{\partial }{\partial x},x\frac{\partial }{\partial y},y\frac{\partial }{\partial y},z\frac{\partial }{\partial y},\frac{\partial }{\partial y},x\frac{\partial }{\partial z},y\frac{\partial }{\partial z},z\frac{\partial }{\partial z},\frac{\partial }{\partial z}\right\}$  \\
 	$(m,n,p)\cdot(x,y,z)=(e^mx,e^my,z+p)$ &$(m,n,p)\cdot(x,y,z)=(e^mx,y+n,z+p)$ \\
$\left\{x\frac{\partial }{\partial x},y\frac{\partial }{\partial x},x\frac{\partial }{\partial y},y\frac{\partial }{\partial y},\frac{\partial }{\partial z}\right\}$ &$\left\{x\frac{\partial }{\partial x},\frac{\partial }{\partial y},\frac{\partial }{\partial z}\right\}$ \\ 	
$\it{Aff}(G,\overline{\nabla}_7)_0\cong \left\{\left[\begin{smallmatrix}
	a&b&0&0\\c&d&0&0\\0&0&1&f
	\end{smallmatrix}\right]\ \bigg|\ ad-bc\ne0  \right\}/\mathbb{Z}^2.$ & $\it{Aff}(\mathbb{T}_8,\overline{\nabla}_8)_0\cong \left\{\left[\begin{smallmatrix}
		a&0&0&0\\0&1&0&b\\0&0&1&c
	\end{smallmatrix}\right]\ \bigg|\ a\ne0 \right\}/\mathbb{Z}^3.$ \\ \hline 
$D_9(a,b,c)=(be^a, e^a,c)$&$D_{10}(a,b,c)=(a+\frac{1}{2}b^2,b, e^c)$ \\
$\left\{x\frac{\partial }{\partial x},y\frac{\partial }{\partial x},z\frac{\partial }{\partial x},\frac{\partial }{\partial x},y\frac{\partial }{\partial y},x\frac{\partial }{\partial z},y\frac{\partial }{\partial z},z\frac{\partial }{\partial z},\frac{\partial }{\partial z}\right\}$ & $\left\{x\frac{\partial }{\partial x},y\dfrac{\partial }{\partial x},z\frac{\partial }{\partial x},\frac{\partial }{\partial x},x\frac{\partial }{\partial y},y\frac{\partial }{\partial y},z\frac{\partial }{\partial y},\frac{\partial }{\partial y},z\frac{\partial }{\partial z}\right\}$\\
$(m,n,p)\cdot(x,y,z)=(e^mx+ne^my,e^my,z+p)$ &$(m,n,p)\cdot(x,y,z)=\left(x+ny+m+\frac{1}{2}n^2,y+n,e^pz\right)$ \\
$\left\{x\frac{\partial }{\partial x}+y\frac{\partial }{\partial y},y\frac{\partial }{\partial x},\frac{\partial }{\partial z}\right\}$ &$\left\{y\frac{\partial }{\partial x}+\frac{\partial }{\partial y},\frac{\partial }{\partial x},z\frac{\partial }{\partial z}\right\}$ \\
$\it{Aff}(\mathbb{T}_9,\overline{\nabla}_9)_0\cong \left\{\left[\begin{smallmatrix}
		a&b&0&0\\0&a&0&0\\0&0&1&c
	\end{smallmatrix}\right]\ \bigg|\ a\ne0  \right\}/\mathbb{Z}^3.$ &$\it{Aff}(\mathbb{T}_{10},\overline{\nabla}_{10})_0\cong  \left\{\left[\begin{smallmatrix}
		1&a&0&b\\0&1&0&a\\0&0&c&0
	\end{smallmatrix}\right]\ \bigg|\ c\ne0  \right\}/\mathbb{Z}^3.$ \\ \hline
$D_{11}(a,b,c)=(a+\frac{1}{2}(b^2+c^2),b, c)$ &$D_{12}(a,b,c)=(a+\frac{1}{2}(b^2-c^2),b, c)$ \\ 
$\it{aff}(\mathbb{R}^3)$&$\it{aff}(\mathbb{R}^3)$\\
$(m,n,p)\cdot(x,y,z)=$ &$(m,n,p)\cdot(x,y,z)=$  \\
	 $\left(x+ny+pz+m+\frac{1}{2}(n^2+p^2),y+n,z+p\right)$ &$\left(x+ny-pz+m+\frac{1}{2}(n^2-p^2),y+n,z+p\right)$\\
$\left\{y\frac{\partial }{\partial x}+\frac{\partial }{\partial y},z\frac{\partial }{\partial x}+\frac{\partial }{\partial z},\frac{\partial }{\partial x}\right\}$ &$\left\{y\frac{\partial }{\partial x}+\frac{\partial }{\partial y},z\frac{\partial }{\partial x}-\frac{\partial }{\partial z},\frac{\partial }{\partial x}\right\}$ \\ 
$\it{Aff}(\mathbb{T}_{11},\overline{\nabla}_{11})_0\cong  \left\{\left[\begin{smallmatrix}
		1&a&b&c\\0&1&0&a\\0&0&1&b
	\end{smallmatrix}\right]\ \bigg|\ a,b,c\in\mathbb{R} \right\}/\mathbb{Z}^3.$ &$\it{Aff}(\mathbb{T}_{12},\overline{\nabla}_{12})_0\cong \left\{\left[\begin{smallmatrix}
		1&a&b&c\\0&1&0&a\\0&0&1&-b
	\end{smallmatrix}\right]\ \bigg|\ a,b,c\in\mathbb{R} \right\}/\mathbb{Z}^3.$\\ \hline 
$D_{13}(a,b,c)=(a+\frac{1}{2}b^2,b, c)$& $D_{14}(a,b,c)=(a+bc+\frac{1}{6}c^3,b+\frac{1}{2}c^2, c)$ \\ 
$\it{aff}(\mathbb{R}^3)$ & $\it{aff}(\mathbb{R}^3)$ \\
$(m,n,p)\cdot(x,y,z)=$
 &$(m,n,p)\cdot(x,y,z)=(x+py+\left(n+\frac{1}{2}p^2\right)z+$ \\ $\left(x+ny+m+\frac{1}{2}n^2,y+n,z+p\right)$ & $+m+np+\frac{1}{6}p^3,y+pz+n+\frac{1}{2}p^2,z+p)$ \\
$\left\{y\frac{\partial }{\partial x}+\frac{\partial }{\partial y},\frac{\partial }{\partial x},\frac{\partial }{\partial z}\right\}$ &$\left\{y\frac{\partial }{\partial x}+z\frac{\partial }{\partial y}+\frac{\partial }{\partial z},z\frac{\partial }{\partial x}+\frac{\partial }{\partial y},\frac{\partial }{\partial x}\right\}$ \\
$\it{Aff}(\mathbb{T}_{13},\overline{\nabla}_{13})_0\cong  \left\{\left[\begin{smallmatrix}
		1&a&0&b\\0&1&0&a\\0&0&1&c
	\end{smallmatrix}\right]\ \bigg|\ a,b,c\in\mathbb{R} \right\}/\mathbb{Z}^3.$ & $\it{Aff}(\mathbb{T}_{14},\overline{\nabla}_{14})_0\cong  \left\{\left[\begin{smallmatrix}
		1&a&b&c\\0&1&a&b\\0&0&1&a
	\end{smallmatrix}\right]\ \bigg|\ a,b,c\in\mathbb{R} \right\}/\mathbb{Z}^3.$ \\ \hline 
\end{tabular}
\vskip5pt 
Finally, the usual torus determined by the action $(m,n,p)\cdot(x,y,z)=(x+m,y+n,z+p)$ of $\mathbb{Z}^3$ on $\mathbb{R}^3$. An easy calculation shows that the subspace of $\it{aff}(\mathbb{R}^3)$ commuting with this action has linear basis given by $\left\{\dfrac{\partial }{\partial x},\dfrac{\partial }{\partial y},\dfrac{\partial }{\partial z}\right\}$, so the Lie group $\it{Aff}(\mathbb{R},\overline{\nabla}_{15}) $ is locally isomorphic to the group  $ \{T:\mathbb{R}^3\rightarrow\mathbb{R}^3\mid T(x,y,z)=(x+a,y+b,z+c)  \}/\mathbb{Z}^3,$ that is, isomorphic to an usual torus. 
More generally, the action $(m_1,\dots,m_n)\cdot(x_1,\dots,x_n)=(x_1+m_1,\dots,x_n+m_n)$   of $\mathbb{Z}^n$ on $\mathbb{R}^n$ determines an n-torus with space of infinitesimal affine transformations isomorphic to the Lie algebra with linear basis $\left\{\dfrac{\partial }{\partial x_1},\dots,\dfrac{\partial }{\partial x_n}\right\}$ and its group of affine transformations $\it{Aff}(\mathbb{R},\overline{\nabla})$  locally isomorphic to the group $ \{T:\mathbb{R}^n\rightarrow\mathbb{R}^n\mid T(x_1,\dots,x_n)=(x_1+a_1,\dots,x_n+a_n),\ a_1,\dots,a_n\in\mathbb{R} \} ,$ also isomorphic to an $n$-dimensional torus.

\section{Further consequences of the results}

We have the following consequences of the results in Section \ref{S:transformationsonhomogeneousspacesfromtransformationsong}
\begin{corollary} \label{existenciarho1rho2} Under the hypothesis of Theorem \ref{C:affinetransformationgrouponhomogeneousspaceswithdiscretequotient}, if $\nabla$ is flat affine and left invariant  there exists a representation $\rho:\widehat{\textit{Aff}}({G/H},\overline{\nabla})\rightarrow\textit{Aff}\left(\mathbb{R}^n,\nabla^0\right)$ with a point of open orbit, where $n=\dim G$, $\nabla^0$ is the usual connection of $\mathbb{R}^n$, and  $\widehat{\textit{Aff}}(G/H,\overline{\nabla})_0$ denote the universal cover of  ${\textit{Aff}}(G/H,\overline{\nabla})_0$.  Furthermore, if the dimension of $V=\it{aff}(G,\nabla^+)^\pi$, the space of projectable complete infinitesimal affine transformations of $G$ relative to $\nabla$, is equal to $n$, then $\rho$ is \'etale.
\end{corollary}
\begin{proof} As $\nabla$ is flat affine and $\pi$ is affine, we have that $\overline{\nabla}$ is flat affine. Corollary 2.4 in \cite{MSV} guarantees the existence of $\rho$ and since right invariant vector fields are projectable, we have that  $\dim V\geq \dim G$, hence Theorem 2.6 in \cite{MSV} guarantees the existence of a point of open orbit. The last statement  follows from the Corollary \ref{C:affinetransformationsofG/HareaffinetransformationsofGcommutingwithRh} and Theorem 2.6 in \cite{MSV}.
\end{proof}
We get a similar result under the hypothesis of Theorem \ref{T:affinetransformationgrouponhomogeneousspaceswithdiscretequotient}, as stated next
\begin{corollary} Under the hypothesis of Theorem \ref{T:affinetransformationgrouponhomogeneousspaceswithdiscretequotient} if $\nabla$ is flat affine there exists a representation $\rho:\widehat{\textit{Aff}}({G/H},\overline{\nabla})\rightarrow\textit{Aff}\left(\mathbb{R}^n,\nabla^0\right)$ with a point of open orbit, where $n=\dim G$.  Furthermore, $\dim V=n$, then $\rho$ is \'etale.
\end{corollary}
 
\begin{corollary}  Under the conditions of Theorem \ref{C:affinetransformationgrouponhomogeneousspaceswithdiscretequotient}, if $(G,\nabla^+)$ is a complete flat affine Lie group then
	\[ \it{aff}(G/H,\overline{\nabla})\cong\{ X\in \mathfrak{a}(G,\nabla^+)\mid X\text{ is }R_h-\text{invariant for all }h\in H \}.  \]
\end{corollary}

\begin{corollary} \label{C:casogrupodelieafinplano} Let $(G,\nabla^+)$ be a  flat affine Lie group with $G$ simply connected and developing map $D:G\rightarrow\mathcal{O}$, then $\it{Aff}(\mathcal{O},\nabla^0)$ is locally isomorphic to the group of transformations of $\it{Aff}(G,\nabla^+)/Aut(D)$ commuting with the right action of $Is_0$ over $G$, where $Is_0$ is the isotropy group at $0$. Moreover, $\it{aff}(G,\nabla^+)^\pi\cong\it{aff}(\mathcal{O},\nabla^0)$
\end{corollary}
\begin{proof} Since $D$ is an affine covering map, the result immediately  follows from Corollary \ref{C:affinetransformationsofG/HareaffinetransformationsofGcommutingwithRh}.  
\end{proof} 
Although the results in this work help to find the group of affine transformations (respectively, the complete infinitesimal affine transformations) of the homogeneous space from the group of affine transformations (respectively, the complete infinitesimal affine transformations) of the Lie group, this corollary shows that sometimes maybe easier to calculate the affine transformations of the homogeneous space as we show in the following example.

\begin{example} Consider $G=\mathbb{R}^2$ with the flat affine left invariant connection $\nabla^+$ determined by affine \'etale representation $\rho:\mathbb{R}^2\longrightarrow\it{Aff}(\mathbb{R}^2)$ with $\rho(x,y)=e^x\begin{bmatrix}	\cos( y) & -sin( y)  \\
		\sin( y) & \ \ cos( y)  \\
	\end{bmatrix}$. The orbit $\mathbb{R}^2\setminus\{(0,0) \}$ is isomorphic to $\mathbb{R}^2/\mathbb{Z}$ where $\mathbb{Z}$ is identified with $H=\{(0,2\pi n)\mid n\in\mathbb{Z} \}$. Using Lemma 2 in \cite{SF}, we get that the group of affine transformations preserving the orbit is $GL_2(\mathbb{R})$ so in local coordinates so, by taking the linear basis (in local coordinates $\{(u,v)\}$) $$\left\{u\dfrac{\partial}{\partial u}+v\dfrac{\partial}{\partial v}, -v\dfrac{\partial}{\partial u}+u\dfrac{\partial}{\partial v},u\dfrac{\partial}{\partial u}-v\dfrac{\partial}{\partial v},v\dfrac{\partial}{\partial u}+u\dfrac{\partial}{\partial v}\right\}$$ of $gl_2(\mathbb{R})\cong \it{aff}(\mathcal{O},\nabla^0)$, an easy calculation shows that the set 
	\[\left\{ \frac{\partial}{\partial x},\dfrac{\partial}{\partial y},
	\sin(2 y)\dfrac{\partial}{\partial x}+\cos(2 y)\dfrac{\partial}{\partial y},
	\cos(2 y)\dfrac{\partial}{\partial x}-\sin(2 y)\dfrac{\partial}{\partial y} \right\}\]is a linear basis for $\it{aff}(G,\nabla^+)^\pi$ in the local coordinates $\{(x,y)\}$.
\end{example}

\end{document}